\documentclass[letterpaper, 10 pt, conference]{ieeeconf}
\IEEEoverridecommandlockouts
\usepackage{graphicx}
\usepackage{subcaption}
\usepackage{amsmath}
\usepackage{mathtools}
\usepackage{dsfont}
\usepackage[utf8]{inputenc}
\usepackage{hyperref}
\usepackage{flushend}

\title{\LARGE \bf
Guaranteed Cost Model Predictive Control-based Driver Assistance System for Vehicle Stabilization Under Tire Parameters Uncertainties}

\author{Carlos M. Massera$^{1}$, Marco H. Terra$^{2}$ and Denis F. Wolf$^{1}$
\thanks{This work was supported by FAPESP grant 2013/24542-7.}
\thanks{$^{1}$Institute of Mathmatics and Computer Science, University of São Paulo, Avenida Trabalhador São-carlense, 400, São Carlos, Brazil {\tt\small massera,denis@icmc.usp.br}}%
\thanks{$^{2}$São Carlos School of Engineering, University of São Paulo, Avenida Trabalhador São-carlense, 400, São Carlos, Brazil {\tt\small terra@sc.usp.br}}%
}

\begin{document}
\maketitle
\thispagestyle{empty}
\pagestyle{empty}

\begin{abstract}
Road traffic crashes have been the leading cause of death among young people. Most of these accidents occur when the driver becomes distracted and a loss-of-control situation occurs. Steer-by-Wire systems were recently proposed as an alternative to mitigate such accidents. This technology enables the decoupling of the front wheel steering angles from the driver hand wheel angle and, consequently, the measurement of road/tire friction limits and the development of novel control systems capable of ensuring vehicle stabilization and safety. However, vehicle safety boundaries are highly dependent on tire characteristics which vary significantly with temperature, wear and the tire manufacturing process. Therefore, design of autonomous vehicle and driver assistance controllers cannot assume that these characteristics are constant or known. Thus, this paper proposes a Guaranteed Cost Model Predictive Controller Driver Assistance System able to avoid front and rear tire saturation and to track the drivers intent up to the limits of handling for a vehicle with uncertain tire parameters. Simulation results show the performance of the proposed approach under time-varying uniformly distributed disturbances.
\end{abstract}

\section{INTRODUCTION}

Road traffic crashes are the leading cause of death among young people between 10 and 24 years old \cite{world2007youth}. 
Most of these accidents occur when the driver is unable to maintain the vehicle control due to fatigue or external factors resulting in loss-of-control scenarios \cite{national2008national}. 
In recent years, both academia and industry have been devoted towards the development of safety systems in order to decrease the number of road accidents. 
% Electronic Stability Control (ESC) is an example of these systems, it reduced the number of fatal crashes in single-vehicle accidents in 32\% \cite{koibuchi1996vehicle}. 
% This system corrects the vehicle yaw rate avoiding under-steer situations, in which the front tire saturates, or over-steer situations, in which the rear tires saturates.

Steer-by-Wire systems were recently introduced in commercial vehicles \cite{ulrich2014top}. 
This system eliminates mechanical coupling between the drivers steering wheel and the front road wheels. 
In 2010, Hsu et al. \cite{hsu2010estimation} demonstrated that a Steer-by-Wire system allowed friction estimation based on the steering torque and Hamann et al. \cite{hamann2014tire} proposed a faster converging method based on the Unscented Kalman Filter. 
The capability to decouple the front wheel steering angles from the driver and to measure saturation limits enable the development of controllers able to predict and avoid vehicle handling saturation limits. 
% Therefore, such controllers could extend the capabilities of current ESC systems and have the potential to further reduce both single-vehicle and multi-vehicle crashes.

Model Predictive Control (MPC) is a class of optimization-based control algorithms which uses an explicit model of the controlled system to predict its future states \cite{badgwell2015model}. 
%Research and commercial interest and usage of such technique increased in recent years, mostly due to its ability to regulate the system state while ensuring both state and input feasibility. This technique has been applied in several different fields, such as: refineries, food processing plants, mining, aerospace and automotive control \cite{qin2003survey}.
% 
Beal and Gerdes \cite{beal2013mpc} proposed a MPC-based Driver Assistance System (DAS) for Steer-by-Wire vehicles able to ensure vehicle handling limits in coordination with a human driver. A linear bicycle model variation, named Affine Force Input (AFI), was also proposed. This model is able to represent the tire dynamics close to a linearization point and to maintain convexity of the optimization problem.
% Carvalho et al. \cite{carvalho2013predictive} designed a MPC for evasive maneuvers based on the iterative linearization of a nonlinear Ackerman model with convexified constraints. 
Bernardini et al. \cite{bernardini2009drive} developed a Hybrid MPC to coordinate steering correction and differential braking. Tire forces were approximated by a piece-wise affine function and the optimization problem of the MPC controller was formulated as a mixed-integer quadratic problem. 
Massera and Denis \cite{massera2015driver} previously proposed a MPC-based DAS for Steer-by-Wire vehicles which ensured vehicle handling limits by incorporating both longitudinal and lateral dynamics by iterative linearizations of a force input nonlinear bicycle model. 

Tire characteristics vary significantly with temperature \cite{tonuk2001prediction}, wear \cite{braghin2006tyre} and their manufacturing process. 
Therefore, the design of autonomous vehicles and driver assistance controllers cannot assume that these characteristics are constant or known. 
The disregard of such uncertainties can lead to poor closed-loop performance of MPCs and, consequently, to the violation of state and control input constraints \cite{rawlings1994nonlinear}. 
Massera et al. \cite{massera2016guaranteed} have proposed a Robust MPC technique, entitled Guaranteed Cost Model Predictive Control (GCMPC), able to guarantee robust stability, robust feasibility and an upper bound to a MPC optimization problem cost for linear system with multiplicative parametric uncertainties.

This paper proposes a Guaranteed Cost Model Predictive Controller Driver Assistance System able to avoid front and rear tire saturation and to track the drivers intent up to the limits of handling for a vehicle with uncertain tire parameters. 
The control law is designed such that both the rear and the front tire never saturate and the vehicle never exceeds its maximum stable yaw rate boundaries.

The remainder of the paper is organized as follows: Section \ref{sec_tire} discusses the tire model and its uncertainties; Section \ref{sec_vehicle} presents the vehicle modeling; Section \ref{sec_controller} describes the proposed controller; Section \ref{sec_experiments} reports the experiments performed; and Section \ref{sec_conclusion} provides the final remarks.

\section{TIRE MODEL}
\label{sec_tire}

There are three categories of models capable of representing the tire dynamics on saturation situations: Finite element analysis \cite{helnwein1993new} \cite{koishi1998tire}; Empirical data approximation, such as the "Magic Tire Formula" \cite{pacejka1992magic}; and the dynamical approximation of the tire by a "brush" model, first proposed by Fiala \cite{fiala1987kraftfahrzeugtechnik}. 
% The third technique, shown in Figure \ref{fig_fiala_model}, provides a good compromise between its ability to describe tires physical properties and its complexity, since it assumes that tires are always at steady state and tire transients are faster than chassis transients. While a linear model provides a good approximation of tire forces for low slip conditions.
The third technique, provides a good compromise between its ability to describe tires physical properties and its complexity, since it assumes that tires are always at steady state and tire transients are faster than chassis transients. While a linear model provides a good approximation of tire forces for low slip conditions.

% \begin{figure}[!ht]
% \centering
% \includegraphics[width=\columnwidth]{images/fiala_model_plot.png}
% \caption{Fiala tire model \cite{fiala1987kraftfahrzeugtechnik} lateral force, as parameterized by \cite{beal2013mpc}, with $ C_{i} = 10^6 $, $ R_{\mu,i} = 0.8 $ and $ \mu = 1 $.}
% \label{fig_fiala_model}
% \end{figure}

\subsection{Linear Tire Model}

A linear approximation of the lateral tire forces is valid for low slip angle situations and enables the development of a linear model of the vehicle dynamics. For a given tire $ i \in \{f, r\} $, the tire lateral force is described by%
\begin{equation}
F_{yi} = - C_i \alpha_i
\label{eq_linear_tire}
\end{equation}%
where $ \alpha_i $ is the tire slip angle and $ C_i $ is the cornering stiffness.

\subsection{Fiala Tire Model}

The Fiala Tire model is parameterized by its slip angle ($ \alpha_{i} $), cornering stiffness ($ C_{i} $), static friction coefficient ($ \mu $), tire normal force ($ F_{zi} $), and the ratio between dynamical and static friction coefficients ($ R_{\mu, i} $) for a given tire $ i \in \{f, r\} $. The lateral tire force for this model is described by%
\begin{multline}
F_{yi} = - f(\alpha_{i}) + \frac{2 - R_{\mu, i}}{3 \mu F_{zi}} | f(\alpha_{i}) | f(\alpha_{i}) - \\ - \frac{1 - \frac{2}{3} R_{\mu, i}}{(3 \mu F_{zi})^{2}} f(\alpha_{i})^3
\label{eq_fiala_unsaturated}
\end{multline}%
for the unsaturated case $\left(|\alpha_{i}| \le tan^{-1}(3 \mu F_{zi} / C_{i}) \right)$ and%
\begin{equation}
F_{yi} = - sign(\alpha_{i}) \mu R_{\mu,i} F_{zi}
\end{equation}%
for the saturated case $\left(|\alpha_{i}| > tan^{-1}(3 \mu F_{zi} / C_{i}) \right)$ where $ f(\alpha_{i}) = C_{i} tan(\alpha_{i}) $. From \eqref{eq_fiala_unsaturated}, the peak lateral force and its slip angle are%
\begin{align}
q &= \left( 1 - \frac{2}{3} R_{mu} \right)^{-1},\\
F_{yi}^{peak} &= \mu F_{zi} \left(- q + \frac{2 - R_{mu}}{3} q^{2} + \frac{1 - \frac{2}{3} R_{mu}}{9} q^{3}\right),\\
\alpha_{i}^{peak} &= tan^{-1}\left( \frac{q \mu F_{zi}}{C_{i}} \right).
\end{align}

\subsection{Parameter Sensitivity}

The normal load $ F_{zi} $ does not vary significantly and it is always coupled with the friction coefficient $ \mu $. Thus, this study will focus on the sensitivity of $ C_i $, $ R_{\mu,i} $ and $ \mu $. 

\begin{figure}[!ht]
\centering
\includegraphics[width=\columnwidth]{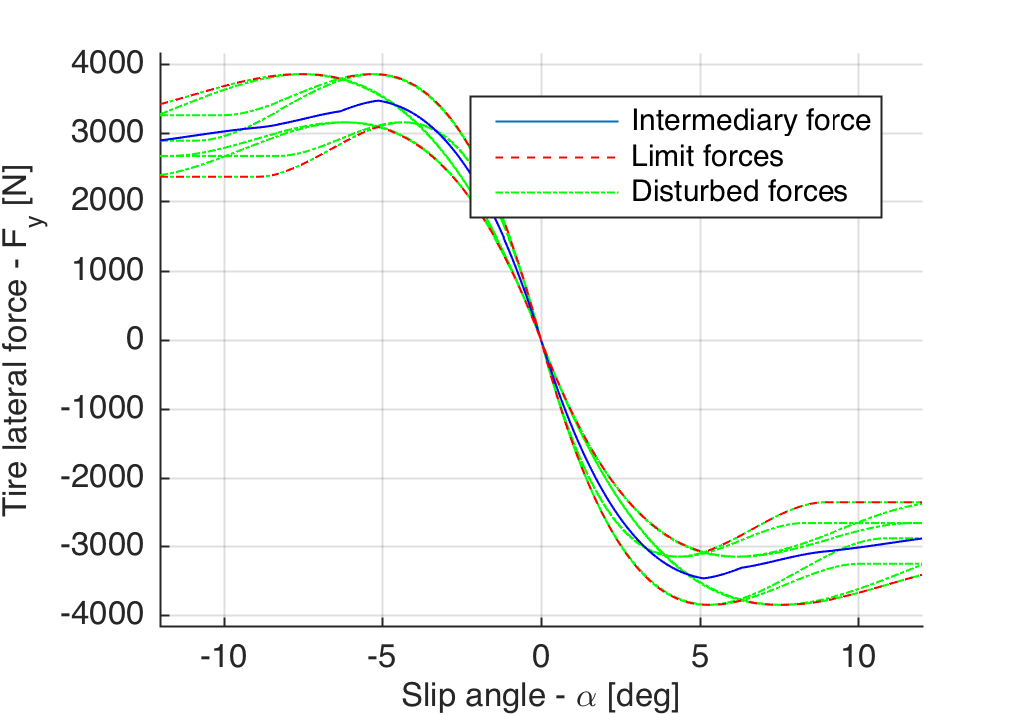}
\caption{Fiala tire model lateral force subject to parametric disturbances of $30\%$ on $C_{i}$, $10\%$ on $R_{\mu,i}$ and $10\%$ on $\mu$.}
\label{fig_fiala_model_sensitivity}
\end{figure}

Let $ p_i = [C_i, R_{\mu,i}, \mu]^T $ be the parameter vector of the $i$-th tire and $ p_i^n $ be the nominal or estimated parameters. Then, the set of possible parametric tire disturbances%
\begin{multline}
\mathcal{W}_i = \{p \mid p \in \mathds{R}^3, C_i \in [0.7, 1.3] C_i^a, \\ R_{\mu,i} \in [0.9, 1.1] R_{\mu,i}^a, \mu \in [0.9, 1.1] \mu^a \}
\end{multline}%
represents uncertainties of $ 30\% $ on $ C_i $, $ 10\% $ on $ R_{\mu,i} $ and $ 10\% $ on $ \mu $. Figure \ref{fig_fiala_model_sensitivity} shows the lateral force profile for each vertex of $ \mathcal{W}_i $, and also presents the upper bound, the lower bound and the intermediary force profiles.

Let $ F_{yi}^{sup}(\alpha_i) = \sup_{p_i \in \mathcal{W}_i} F_{yi}(\alpha_i) $ and $ F_{yi}^{inf}(\alpha_i) = \inf_{p_i \in \mathcal{W}_i} F_{yi}(\alpha_i) $ be the upper and lower bound forces for a given slip angle, respectively. Then, the intermediary force $ \bar{F}_{yi}(\alpha_i) $ and the force deviation are defined by
\begin{align}
\bar{F}_{yi}(\alpha_i) =& \left( F_{yi}^{sup}(\alpha_i) + F_{yi}^{inf}(\alpha_i) \right) / 2 \\
\partial F_{yi}(\alpha_i) =& \left( F_{yi}^{sup}(\alpha_i) - F_{yi}^{inf}(\alpha_i) \right) / 2
\label{eq_intermadiary_force}
\end{align}%
From \eqref{eq_intermadiary_force} we conclude that%
\begin{multline}
\forall p_i \in \mathcal{W}_i \exists \gamma_f \in \mathds{R}, |\gamma_f| \le 1: \\
F_{yi}(\alpha_i) = \bar{F}_{yi}(\alpha_i) + \gamma_f \partial F_{yi}(\alpha_i).
\label{eq_fy_uncertainty}
\end{multline}%
Considering definitions similar to \eqref{eq_intermadiary_force} for the local cornering stiffness $ \widehat{C}_i(\alpha_i) = - \partial F_{yi}(\alpha_i) / \partial \alpha_i $ results%
\begin{multline}
\forall p_i \in \mathcal{W}_i \exists \gamma_c \in \mathds{R}, |\gamma_c| \le 1: \\ \widehat{C}_{i}(\alpha_i) = \bar{C}_{i}(\alpha_i) + \gamma_c \partial C_{i}(\alpha_i).
\label{eq_c_uncertainty}
\end{multline}

\section{VEHICLE MODEL}
\label{sec_vehicle}

\begin{figure}[!ht]
\centering
\includegraphics[width=0.8\columnwidth]{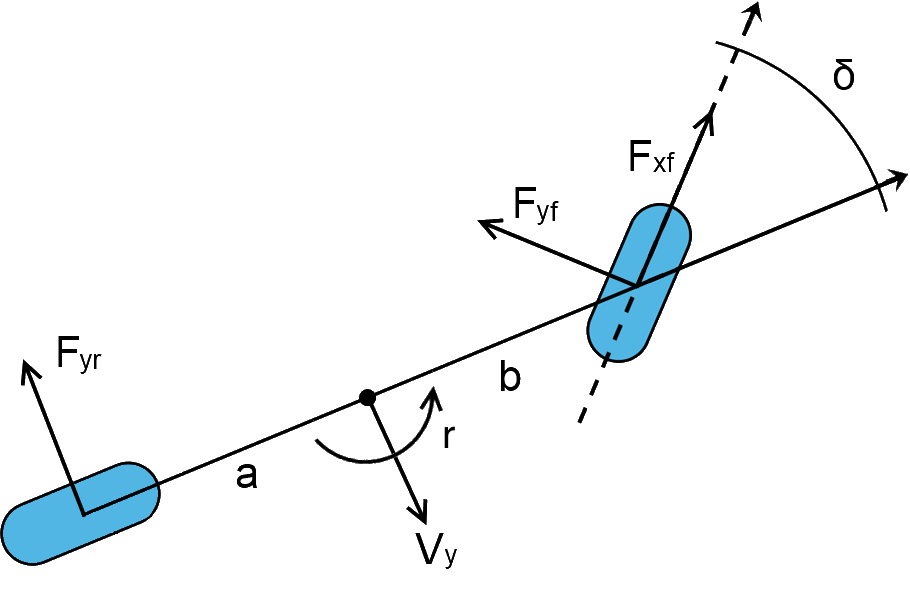}
\caption{Bicycle model representation}
\label{fig_bicycle_model}
\end{figure}

This section describes the linear bicycle model, the AFI bicycle model, proposed by Beal and Gerdes \cite{beal2013mpc}, and the proposed parametric uncertain AFI bicycle model. 
The bicycle model is a simplification of the vehicle dynamics where the left and right front wheels are replaced by a virtual front wheel in the center of the front axle and, analogously, the left and right rear wheels are replaced by a virtual wheel in the center of the rear axle. A representation of a vehicle on this form is shown in Figure \ref{fig_bicycle_model}.

\subsection{Linear Bicycle Model}

The linear dynamical model is a two-state model, which describes both yaw-rate $ r $ and lateral speed $ v_y $ of a vehicle based on the bicycle model [5]. It uses low angle and constant longitudinal speed assumptions ($cos(x) \approx 1$, $sin(x) \approx x$, $tan(x) \approx x$ and $\dot{v}_x \approx 0$) \cite{massera2015driver}. Its equations of motion are given by%
\begin{equation}
\dot{v}_y = \frac{F_{yf} + F{yr}}{m} - v_x r, \; \dot{r} = \frac{a F_{yf} - b F_{yr}}{I_z}
\end{equation}%
in which $ a $, $ b $, $ m $ and $ I_z $ are the distance from the front axle to the center of gravity, the distance from the rear axle to the center of gravity, the vehicle mass and the vehicle yaw moment of inertia, respectively. The tire forces are described by the linear tire model from \eqref{eq_linear_tire} with approximate slip angles%
\begin{equation}
\alpha_f \approx \frac{v_y + a r}{v_x} - \delta, \; \alpha_r \approx \frac{v_y - b r}{v_x}
\end{equation}
in which $ delta $ is the steering angle. Therefore, the linear bicycle system can be represented in state-space form $ \dot{x} = A_l x + B_l u $ with $ x = [v_y, r]^T $, $ u = \delta $ and%
\begin{equation}
A_l = \begin{bmatrix}
- \frac{C_f + C_r}{m v_x} & - \frac{a C_f - b C_r}{m v_x} - v_x \\ 
- \frac{a C_f - b C_r}{I_z v_x} & - \frac{a^2 C_f + b^2 C_r}{I_z v_x}
\end{bmatrix}, \;
B_l = \begin{bmatrix}
\frac{C_f}{m}\\ 
\frac{a C_f}{I_z}
\end{bmatrix}
\end{equation}

\subsection{Affine Force Input Bicycle Model}

The AFI bicycle model was first proposed by Beal and Gerdes \cite{beal2013mpc}. The front tire forces are abstracted from the system dynamics, such that the system input is chosen to be $ F_{yf} $, instead of $ \delta $. Also, the rear tire forces are approximated by a first order Taylor series%
\begin{equation}
F_{yr}(\alpha_r) \approx F_{yr}(\widehat{\alpha}_r) - (\alpha_r - \widehat{\alpha}_r) \widehat{C}_r(\widehat{\alpha}_r).
\end{equation}

Thus, the AFI bicycle system can be represented in the affine state-space form $ \dot{x} = A_{afi} x + B_{afi} u + c_{afi} $ with $ x = [v_y, r]^T $, $ u = F_{yf} $ and%
\begin{multline}
A_{afi} = \begin{bmatrix}
- \frac{\widehat{C}_r(\widehat{\alpha}_r)}{m v_x} & \frac{b \widehat{C}_r(\widehat{\alpha}_r)}{m v_x} - v_x \\ 
- \frac{b \widehat{C}_r(\widehat{\alpha}_r)}{I_z v_x} & - \frac{b^2 \widehat{C}_r(\widehat{\alpha}_r)}{I_z v_x}
\end{bmatrix}, \\
B_{afi} = \begin{bmatrix}
\frac{1}{m}\\ 
\frac{a}{I_z}
\end{bmatrix}, \;
c_{afi} = \begin{bmatrix}
\frac{\widehat{F}_{yr}(\widehat{\alpha}_r) + \widehat{C}_r(\widehat{\alpha}_r) \widehat{\alpha}_r}{m}\\
-b \frac{\widehat{F}_{yr}(\widehat{\alpha}_r) + \widehat{C}_r(\widehat{\alpha}_r) \widehat{\alpha}_r}{I_z}.
\end{bmatrix}
\label{eq_afi_model}
\end{multline}

This model is able to represent the vehicle lateral dynamics up to the limits of handling. However, it is only valid for rear slip angles close to the linearization point $ \hat{\alpha_r} $.

\subsection{Parametric Uncertain AFI Bicycle Model}

Disturbance on tire parameters has direct influence on the system dynamics of the AFI model. Front tire disturbances causes a multiplicative uncertainty on the input $ u $, while rear tire disturbances results in multiplicative uncertainties on matrices $A$ and $c$.

The substitution of the parametrically uncertain tire model from \eqref{eq_fy_uncertainty} and \eqref{eq_c_uncertainty} on the AFI system dynamics from \eqref{eq_afi_model} yields%
\begin{equation}
\begin{matrix}
A_{afi} = \begin{bmatrix}
- \frac{\bar{C}_{i}(\widehat{\alpha}_r) + \gamma_c \partial C_{r}(\widehat{\alpha}_r)}{m v_x} & \frac{b \left[\bar{C}_{r}(\widehat{\alpha}_r) + \gamma_c \partial C_{r}(\widehat{\alpha}_r)\right]}{m v_x} - v_x \\ 
- \frac{b \left[\bar{C}_{r}(\widehat{\alpha}_r) + \gamma_c \partial C_{r}(\widehat{\alpha}_r)\right]}{I_z v_x} & - \frac{b^2 \left[\bar{C}_{r}(\widehat{\alpha}_r) + \gamma_c \partial C_{r}(\widehat{\alpha}_r)\right]}{I_z v_x}
\end{bmatrix}\\
B_{afi} = \left(1 + \gamma_b \frac{\partial F_{yf}(\widehat{\alpha}_f)}{\bar{F}_{yf}(\widehat{\alpha}_f)} \right) \begin{bmatrix}
\frac{1}{m}\\ 
\frac{a}{I_z}
\end{bmatrix}\\
c_{afi} = \begin{bmatrix}
\frac{\bar{F}_{yr}(\widehat{\alpha}_r)) + \gamma_f \partial F_{yr}(\widehat{\alpha}_r) + \left[\bar{C}_{r}(\widehat{\alpha}_r) + \gamma_c \partial C_{r}(\widehat{\alpha}_r)\right] \widehat{\alpha}_r}{m}\\
-b \frac{\bar{F}_{yr}(\widehat{\alpha}_r)) + \gamma_f \partial F_{yr}(\widehat{\alpha}_r) + \left[\bar{C}_{r}(\widehat{\alpha}_r) + \gamma_c \partial C_{r}(\widehat{\alpha}_r)\right] \widehat{\alpha}_r}{I_z}.
\end{bmatrix}
\end{matrix}
\label{eq_pu_afi_model}
\end{equation}

This system dynamics can be rewritten in the parametric uncertain affine state-space form%
\begin{multline}
\dot{x} = (\bar{A}_{afi} + H_{afi} \Delta_{afi} E_{a,afi}) x + (\bar{B}_{afi} + H_{afi} \Delta E_{b,afi}) u + \\ + (\bar{c}_{afi} + H_{afi} \Delta_{afi} E_{c,afi})
\end{multline}%
with%
\begin{equation}
\bar{A}_{afi} = \begin{bmatrix}
- \frac{\bar{C}_{i}(\widehat{\alpha}_r)}{m v_x} & \frac{b \bar{C}_{i}(\widehat{\alpha}_r)}{m v_x} - v_x \\ 
- \frac{b \bar{C}_{i}(\widehat{\alpha}_r)}{I_z v_x} & - \frac{b^2 \bar{C}_{i}(\widehat{\alpha}_r)}{I_z v_x}
\end{bmatrix}
\end{equation}
\begin{equation}
\bar{B}_{afi} = \begin{bmatrix}
\frac{1}{m}\\ 
\frac{a}{I_z}
\end{bmatrix}
\end{equation}
\begin{equation}
\bar{c}_{afi} = \begin{bmatrix}
\frac{\bar{F}_{yr}(\widehat{\alpha}_r) + \bar{C}_r(\widehat{\alpha}_r) \widehat{\alpha}_r}{m}\\
-b \frac{\bar{F}_{yr}(\widehat{\alpha}_r) + \bar{C}_r(\widehat{\alpha}_r) \widehat{\alpha}_r}{I_z}
\end{bmatrix}
\end{equation}
\begin{equation}
H_{afi} = \begin{bmatrix}
m^{-1} & m^{-1} \\
- b I_z^{-1} & a I_z^{-1}
\end{bmatrix}
\end{equation}
\begin{equation}
E_{a,afi} = \begin{bmatrix}
- v_x^{-1} \partial C_r(\widehat{\alpha}_r) & b v_x^{-1} \partial C_r(\widehat{\alpha}_r) \\
0 & 0 \\
0 & 0
\end{bmatrix}
\end{equation}
\begin{equation}
E_{b,afi} = \begin{bmatrix}
0 & 0 & \partial F_{yf}
\end{bmatrix}^T
\end{equation}
\begin{equation}
E_{c,afi} = \begin{bmatrix}
\partial C_r(\widehat{\alpha}_r) \widehat{\alpha}_r & \partial F_{yr}(\widehat{\alpha}_r) & 0
\end{bmatrix}^T
\end{equation}
\begin{equation}
\Delta_{afi} = \begin{bmatrix}
\gamma_f & \gamma_c & 0 \\ 
0 & 0 & \gamma_g
\end{bmatrix}
\end{equation}

\section{CONTROLLER DESIGN}
\label{sec_controller}

The controllers should primarily ensure driver and vehicle safety, avoiding unsafe yaw rates and tire saturation. Secondarily, it should follow the drivers intent as close as possible and avoid uncomfortable and unnecessary interventions.

\subsection{System and Functional Modeling}

% The system model must accurately represent vehicle dynamics up to the handling limits and correctly predict the drivers intent for a short horizon. 
Since the majority of drivers guide their vehicles on state regions where the linear model is valid, the linear bicycle model was chosen to represent the drivers intent for the vehicle. 
Whereas the parametric uncertain affine force input model was chosen to depict the vehicle dynamics, since it correctly models the vehicle dynamics up to the limits of handling and account for tire parameter uncertainties.

The system state was defined as $ x = [v_y^{ref}, r^{ref}, \delta^{ref}, v_y, r, 1]^T $ where $ v_y^{ref} $ and $ r^{ref} $ are the reference lateral velocity and yaw rate, respectively, $ \delta^{ref} $ is the drivers hand wheel angle, $ v_y $ and $ r $ are the effective vehicle lateral velocity and yaw rate, respectively, and there is an one-valued state to incorporate the affine terms of the AFI model. The control input was defined as $ u = F_{yf} $. Therefore, the continuous time system dynamics can be described as $ \dot{x} = (A + H_c \Delta E_a) x + (B + H_c \Delta E_b) u $ with%
\begin{equation}
A = \begin{bmatrix}
A_l & B_l & \multicolumn{2}{c}{\mathds{O}_{2,3}} \\
\multicolumn{1}{c}{\mathds{O}_{1,2}} & 1 & \multicolumn{2}{c}{\mathds{O}_{1,3}} \\
\multicolumn{2}{c}{\mathds{O}_{3,3}} & \bar{A}_{afi} & \bar{c}_{afi} \\
\multicolumn{3}{c}{\mathds{O}_{1,5}} & 1
\end{bmatrix}, \; 
B = \begin{bmatrix}
0 \\
0 \\
\bar{B}_{afi}\\
0
\end{bmatrix},
\label{eq_system_def}
\end{equation}%
$ H_c = [\mathds{O}_{2,3}, H_{afi}^T, \mathds{O}_{2,1}]^T $, $ \Delta = \Delta_{afi} $, $ E_a = [\mathds{O}_{3,3}, E_{a,afi}, E_{c,afi}] $ and $ E_b = E_{b,afi} $ in which $ \mathds{O}_{i,j} $ is a zero-valued matrix with $ i $ rows and $ j $ columns.

This system was discretized with a first-order hold on $ w = \Delta E_a x + \Delta E_b u $ such that%
\begin{equation}
\begin{bmatrix}
F & G & H\\
\mathds{O}_{3,6} & \multicolumn{2}{c}{\mathds{I}_{3,3}}
\end{bmatrix} = exp \left( T_s\begin{bmatrix}
A & B & H_c\\
\multicolumn{3}{c}{\mathds{O}_{3,9}}
\end{bmatrix} \right)
\end{equation}%
in which $T_s$ is the sampling time in seconds and $ \mathds{I}_{i,j} $ is an identity matrix with $ i $ rows and $ j $ columns.

The controller objective was designed as the quadratic minimization of $ v_y^{ref} - v_y $, $ r^{ref} - r $ and $ F_{yf} $ in discrete time. Thus, it was represented in the form%
\begin{equation}
J(x_0, \mathbf{u}, N) = x_N^T S x_N + \underset{k = 0}{\overset{N}{\sum}} x_k^T Q x + u_k^T R u_k
\end{equation}%
where $ S $ is the terminal cost,%
\begin{equation}
Q = \begin{bmatrix}
W_{v_y} & 0 & 0 & -W_{v_y} & 0 & 0\\
0 & W_r & 0 & 0 & - W_r & 0\\
0 & 0 & 0 & 0 & 0 & 0\\
-W_{v_y} & 0 & 0 & W_{v_y} & 0 & 0\\
0 & -W_r & 0 & 0 & W_r & 0\\
0 & 0 & 0 & 0 & 0 & 0\\
\end{bmatrix}
\end{equation}%
and $ R = W_{F_{yf}} $ with $ W_{v_y}, W_r, W_{F_{yf}} > 0 $ such that $ Q \succeq 0 $ and $ R \succ 0 $.

\subsection{Safety Handling Constraints}

Vehicle loss-of-control situations occur when one or more tires saturate, therefore maintaining vehicle controllability is analogous to avoiding tire saturation.

% Beal and Gerdes \cite{beal2013mpc} proposed a Envelope Controller where a constraints to limit the rear slip angle (dashed red lines) and another to limit the yaw rate (blue lines) are used for rear wheel saturation avoidance (shown in Figure \ref{fig_rear_boundaries}).

Beal and Gerdes \cite{beal2013mpc} proposed a Envelope Controller where a constraints to limit the rear slip angle and another to limit the yaw rate are used for rear wheel saturation avoidance.
% 
% \begin{figure}[!ht]
% \centering
% \includegraphics[width=\columnwidth]{images/boundaries.png}
% \caption{Yaw rate and rear slip constraints for $ v_x = 20 m/s $ and $ \mu = 1 $}
% \label{fig_rear_boundaries}
% \end{figure}
% 
From the slip angle definition, the rear slip limit is%
\begin{equation}
\left| tan^{-1} \left( \frac{v_{y} - b r}{v_{x}} \right) \right| \le \alpha_r^{max}
\label{eq_slip_constraint}
\end{equation}%
in which $ \alpha^{max} = \alpha^{peak} $ to avoid skidding scenarios. Therefore, the rear slip constraint is represented by%
\begin{equation}
- v_x tan(\alpha_r^{peak}) \le v_{y} - b r \le v_x tan(\alpha_r^{peak}).
\end{equation}

The yaw rate constraint ensures the invariability of the operation envelope. It guarantees that the yaw rate value can be maintained at steady-state. Thus, the maximum yaw rate is \cite{massera2015driver}%
\begin{equation}
r^{max} = \frac{\mu}{v_{x}}\:\frac{ab + max(a,b)^{2}}{ min(a,b) (a + b)}
\end{equation}%
and the yaw rate constraint is%
\begin{equation}
- r^{max} \le r \le r^{max}.
\end{equation}

The force limits for the front tire are given by%
\begin{equation}
F_{xf}^{2} + F_{yf}^{2} \le \mu ^ 2 F_{zf} ^ 2.
\end{equation}%
Since lateral forces have higher influence on vehicle safety the longitudinal force is assumed to be chosen such that%
\begin{equation}
|F_{xf}| \le \sqrt{\mu ^ 2 F_{zf} ^ 2 - F_{yf}^{2}}.
\end{equation} %
Therefore, the front tire lateral force is constrained by%
\begin{equation}
- \mu F_{zf} \le F_{yf} \le \mu F_{zf}
\end{equation}%
and the longitudinal forces will use any remainder available force.

All constraints presented are represented by their generic form $ M x_k + N u_k \le o $ where%
\begin{equation}
M = \begin{bmatrix}
0 & 0 & 0 & 1 & -b & 0\\
0 & 0 & 0 & -1 & b & 0\\
0 & 0 & 0 & 0 & 1 & 0\\
0 & 0 & 0 & 0 & -1 & 0\\
0 & 0 & 0 & 0 & 0 & 0\\
0 & 0 & 0 & 0 & 0 & 0\\
\end{bmatrix},\;
N = \begin{bmatrix}
0\\
0\\
0\\
0\\
1\\
-1
\end{bmatrix}
\end{equation}%
and%
\begin{equation}
o = \begin{bmatrix}
1 & 1 & 0 & 0 & 0 & 0\\
0 & 0 & 1 & 1 & 0 & 0\\
0 & 0 & 0 & 0 & 1 & 1\\
\end{bmatrix}^T \begin{bmatrix}
v_x tan(\alpha_r^{peak}) \\
r^{max} \\
\mu F_{zf}
\end{bmatrix}.
\end{equation}

\subsection{Guaranteed Cost Model Predictive Controller}

\begin{figure*}[t!]
\centering
\includegraphics[width=\textwidth]{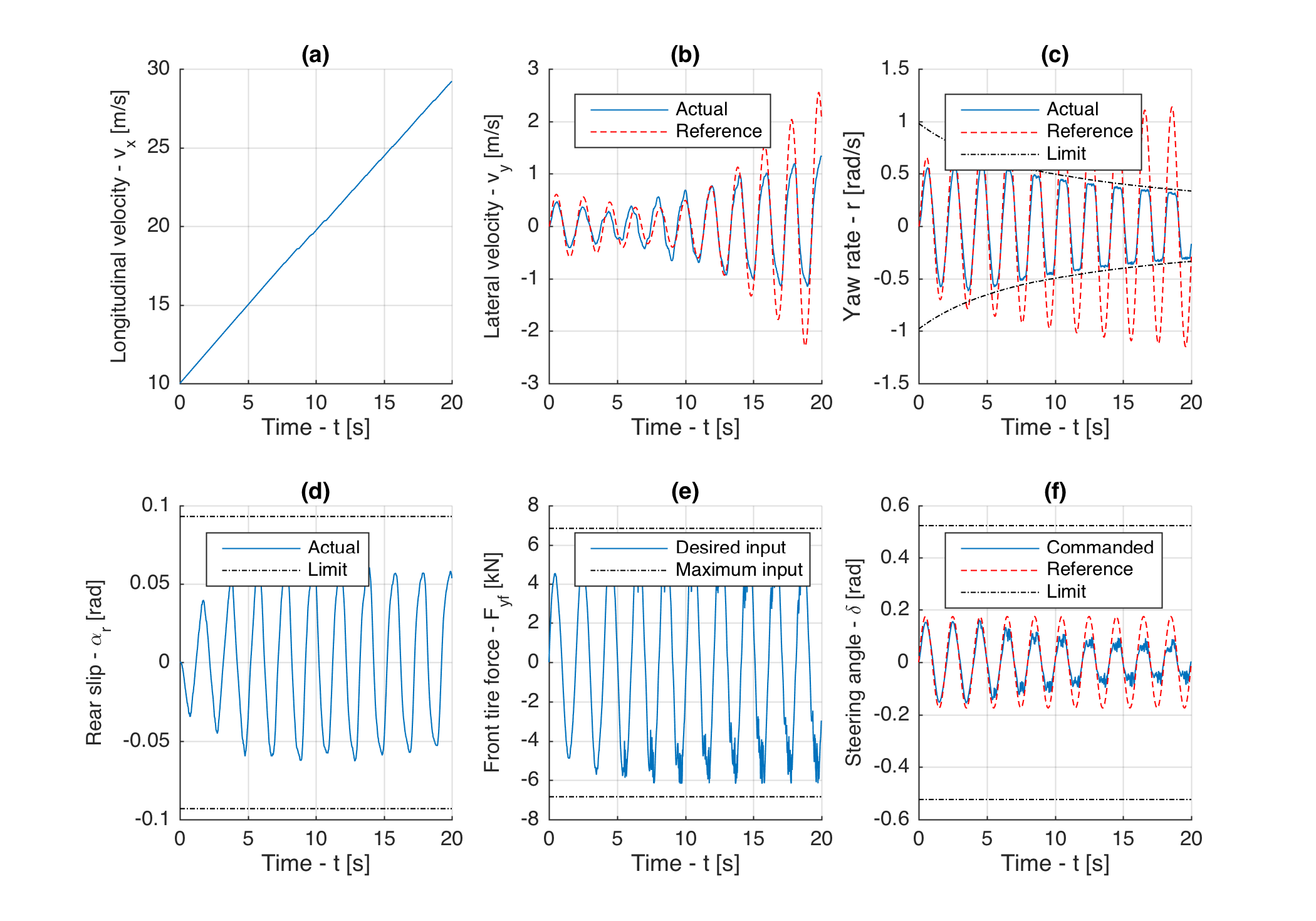}
\caption{Simulation results for a slalom maneuver. (a), (b) and (c) present the longitudinal velocity, lateral velocity and yaw rate, respectively. (d), (e) and (f) show the rear slip angle, the desired front tire force and the commanded steering angle, respectively.}
\label{fig_results}
\end{figure*}

The GCMPC is defined by the optimization problem%
\begin{equation}
\begin{matrix*}[l]
J^*(x_0, N) < & \underset{\mathbf{v}}{\inf} & \underset{k = i}{\overset{N-1}{\sum}} v_k \bar{R}_k v_k + x_i^T S_0 x_i\\
& s.t. & x_{k+1} = (F - G K_k) x_k + G v_k\\
& & (M^{(i)} - N^{(i)} K_k) x_k + N^{(i)} v_k +\\
& & \phantom{(M^{(i)} - N} + \bar{\Phi}_{k,i}(\mathbf{x}, \mathbf{v}) \le o^{(i)}
\end{matrix*}
\label{eq_gcmpc_opt}
\end{equation}%
where $ M^{(i)} $, $ N^{(i)} $, $ o^{(i)} $ are the $i$-th row of matrices $M$, $N$ and vector $o$, respectively, $ S_0 $ is the guaranteed cost matrix, $ K_k $ is the guaranteed cost control at timestep $ k $, $ u_k = - K_k x_k + v_k $, $ \bar{\Phi}_{k,i}(\mathbf{x}, \mathbf{v}) $ is the robustness margin, defined by%
\begin{align}
\bar{\Phi}_{k,i}(\mathbf{x}, \mathbf{v}) &= \underset{j = 0}{\overset{k - 1}{\sum}} || \widetilde{A}_k^{(i)} \widetilde{F}^{k - j - 1} H ||_1 \bar{\phi}_{j,1}(\mathbf{x}, \mathbf{v}),\\
\bar{\phi}_k(\mathbf{x}, \mathbf{v}) &= \phi_k(x_k, v_k) + \underset{i = 0}{\overset{k - 1}{\sum}} c(k,i) \phi_i(x_i, v_i),\\
c(k,i) &= \rho_{k-i-1} + \underset{j = 0}{\overset{k - i - 2}{\sum}} \rho_j c(k-j-1, i), 
\end{align}%
where $\phi_k(x, v) = || \widetilde{E}_{1,k} x + E_2 v ||_2$ is the disturbance boundary, $ \rho_i = || E_1 \widetilde{F}^i H ||_2 $ is the disturbance propagation norm, $ \mathbf{x} = \{x_k \mid k \in [0, N]\} $ is the sequence of states, $ \mathbf{u} = \{u_k \mid k \in [0, N-1]\} $ is the sequence of inputs, $ \mathbf{v} = \{v_k \mid k \in [0, N-1]\} $ is the sequence of feasibility offsets, $ \bar{R}_k = R + \epsilon^{-1} E_b^T E_b + G^T X_{k+1} G $ is the offset cost matrix, $ X_{k} = \left( S^{-1} + \epsilon H H^T \right)^{-1} $ , $ \epsilon > 0 $ is obtained such that $ \text{tr}(S) $ is bounded and minimal and $ J^*(x_0, N) $ is the cost of the parametric uncertain MPC optimization problem%
\begin{equation}
\begin{matrix*}[l]
J^*(x_0, N) = & \underset{\mathbf{u}}{\inf} & J(x_0, \mathbf{u}, N)\\
& s.t. & x_{k+1} = F x_k + G u_k\\
& & M x_k + N u_k \le o
\end{matrix*}
\end{equation}

For further details on the method and its robustness proofs, the reader should refer to \cite{massera2016guaranteed} and to the numerical example implementation available online\footnotemark[1].

\footnotetext[1]{The numerical example source code for the GCMPC is available at: \url{https://github.com/cmasseraf/gcmpc}}

\subsection{Driver Assistance Controller Implementation}

Since the system dynamics depend on $ v_x $ and $ \hat{\alpha}_r $, gain scheduling was used to obtain matrices $ F $, $ G $, $ E_a $, $ S_0 $, $ X $, $ K_k $ and $ \widetilde{K} $. These values were then used in conjunction with the current state $ x_0 $  as inputs to a multi-parametric Quadratic Program (mpQP) designed using YALMIP \cite{lofberg2004yalmip} and Gurobi \cite{gurobi} to represent the optimization problem from \eqref{eq_gcmpc_opt}. The desired front tire force was obtained as $ F_{yf} = - K_0 x_0 + v_0 $ and the steering angle was obtained based on the inverse intermediary front tire lateral force $ \left( \alpha_f = \bar{F}_{yf}^{-1}(F_{yf}) \right) $%
\begin{equation}
\delta = tan^{-1}\left( \frac{v_y + a r }{v_x} \right) - \bar{F}_{yf}^{-1}(F_{yf}).
\end{equation}%
Since it is not desired to saturate the front tire, this inverse function considers only the region where the intermediary tire force is monotonic $ \left(\alpha_f \in \left[- \alpha_f^{peak}, \alpha_f^{peak}\right] \right)$.

Table \ref{tab_controller} presents the parameters values used for the proposed controller.

\begin{table}[thpb]
	\centering
	\caption{Controller Parameters}
	\label{tab_controller}
    \renewcommand{\arraystretch}{1.2}
    \begin{tabular}{|lccc|}    
    \hline
    Parameter  & Symbol  & Value  & Unit \\
    \hline
    Prediction horizon & $ N $ & $ 15 $ & $ - $ \\
    Lateral Speed weight & $ W_{v_{y}} $ & $ 1 $ & $ s^{2} / m^{2} $ \\
    Yaw rate weight & $ W_{r} $ & $ 10^6 $ & $ s^{2} / rad^{2} $ \\
    Lateral force weight & $ W_{F_{yf}} $ & $ 10^{-10} $ & $ 1 / N^{2} $\\
    Maximum yaw rate & $ r^{max} $ & \footnotemark[2] & $ rad / s $ \\
    Maximum rear slip angle & $ \alpha_r^{peak} $ & \footnotemark[2] & $ rad $ \\
    Sampling Time & $T_s$ & $0.02$ & s\\
    \hline
    \end{tabular}
\end{table}%
\footnotetext[2]{Parameter dynamically evaluated during control loop}

\section{SIMULATION RESULTS}
\label{sec_experiments}

The simulation scenario used for validation consists of a vehicle at $ 10 m/s $ accelerating at $ 1 m /s^2 $ for $ 20 s $ on a slalom maneuver, where the hand wheel angle is a sine with $ 0.5 Hz $ frequency and $ 10 deg $ amplitude. Table \ref{tab_vehicle_parameters} presents the simulated vehicle parameters.

\begin{table}[thpb]
	\centering
	\caption{Simulated vehicle parameters}
	\label{tab_vehicle_parameters}
    \renewcommand{\arraystretch}{1.2}
    \begin{tabular}{|lccc|}    
    \hline
    Parameter  & Symbol  & Value & Unit \\
    \hline
    Vehicle mass & $ m $ & 1231 & $ Kg $ \\
    Inertia Moment $ z $ & $ I_{z} $ & 2034.5 & $ \frac{Kg}{m^{2}} $ \\
    Front axle distance to CG & $ a $ & 1.07 & $ m $ \\
    Back axle distance to CG & $ b $ & 1.40 & $ m $ \\
    Front tire cornering stiffness & $ C_{f} $ & 100000 & $ \frac{N}{rad} $ \\
    Rear tire cornering stiffness & $ C_{r} $ & 130000 & $ \frac{N}{rad} $ \\
    Tire-road friction coefficient  & $ \mu $ & 1 & - \\
    Friction coefficient ratio & $R_{\mu,i}$ & 0.8 & - \\
    Cornering stiffness uncertainty & - & 20 & \% \\
    Tire-road friction uncertainty & - & 10 & \% \\
    Friction ratio uncertainty & - & 10 & \% \\
    \hline
    \end{tabular}
\end{table}

Simulation results are shown in Figure \ref{fig_results}, where the tire parameters were chosen to be time-varying and uniformly distributed within the predetermined boundaries. The initial state of the evaluated scenario was $ x = [0, 0, 0, 0, 0, 1] $ and $ v_x = 10 $. 

Since the vehicle has natural under-steer behavior ($ a < b $), the yaw rate limits were reached before rear slip and front tire forces limits. The yaw rate boundary varies significantly since it is inversely dependent on longitudinal velocity. We observed the constraints are satisfied and the vehicle maintained safe operation at all times. Also, a feasibility robustness margin exists even at saturating cases (Fig. \ref{fig_results}c at $ t = 10 s $). It is also possible to notice that tire uncertainties cause higher amplitude disturbances on high slip angles cases rather than on small slip angles, mostly due to the significant variation in the affine term present in the dynamics (Fig. \ref{fig_results}e and \ref{fig_results}f at $ t > 7 s $).

\section{CONCLUSIONS}
\label{sec_conclusion}

This paper proposed a Steer-by-Wire-based lateral dynamics driver assistance system using Guaranteed Cost Model Predictive Control. Front and rear tire saturation limits where defined as operation boundaries in order to avoid loss-of-control situations. Stability and feasibility robustness to tire parameter uncertainties were also incorporated to the system design to ensure safe operation on most driving conditions, since tire parameters vary with temperature, wear and manufacturing process; and cannot be assumed constant or known.

Future work on the controller will focus on the incorporation of static and dynamic objects to the controller safety envelope constraints and the mitigation of tire blowouts.

\bibliographystyle{IEEEtranS}
\bibliography{refs_paper}

\end{document}